\documentclass[12pt]{amsart}
\usepackage{amsmath,amsthm,amssymb, bm}
\usepackage[all]{xy}

\newif\ifappendix
\appendixfalse

\newcommand{\Pa}[9]{\bibitem{#1} {#2}, \emph{#3}, {#4} \textbf{#5} ({#6}), {#7}--{#8}.}
\newcommand{\hl}[1]{{}\save []+<0.5cm,0.15cm>*\txt{#1} \restore}

\newcommand{\sms}{{\ }}
\newcommand{\intvl}[2]{{[#1(#2),#1(#2\!+\!1))}}

\newcommand{\bo}[1]{{\bm{[}\sms{}#1\sms{}\bm{]}}}
\newcommand{\rest}{\upharpoonright}
\newcommand{\alephes}{{\aleph_0}}

\newcommand{\FinSeqs}[1]{{{#1}^{<\alephes}}}
\newcommand{\vphi}{\varphi}
\newcommand{\Bgp}{{\Z^\N}}
\newcommand{\Cgp}{{\Z_2^{\ \N}}}

\newcommand{\arx}[1]{\texttt{http://arxiv.org/abs/#1}}
\newcommand{\ed}{

\end{document}}

\newcommand{\inv}{^{-1}}
\newcommand{\Cantor}{{\{0,1\}^\N}}

\newcommand{\CH}{the Continuum Hypothesis}
\newcommand{\MA}{Martin's Axiom}

\newcommand{\im}{\op{im}}

\newcommand{\Union}{\bigcup}

\newcommand{\fb}{\mathfrak{b}}
\newcommand{\fd}{\mathfrak{d}}
\newcommand{\fc}{\mathfrak{c}}

\newcommand{\seq}[1]{\{#1\}_{n\in\N}}
\newcommand{\setseq}[1]{\{#1 : n\in\N\}}

\newcommand{\cM}{\mathcal{M}}
\newcommand{\cU}{\mathcal{U}}

\newcommand{\mbx}[1]{{\mbox{#1}}}

\long\def\forget#1\forgotten{}

\newcommand{\Fin}{{[\N]^{<\alephes}}}
\newcommand{\PN}{{P(\N)}}
\newcommand{\roth}{{[\N]^{\alephes}}}

\newcommand{\Impl}{\Rightarrow}
\newcommand{\Iff}{\Leftrightarrow}

\newcommand{\cF}{\mathcal{F}}

\newcommand{\N}{\mathbb{N}}
\newcommand{\Z}{\mathbb{Z}}
\newcommand{\NN}{{\N^\N}}

\newcommand{\NNup}{{\N^{\uparrow\N}}}

\newcommand{\R}{\mathbb{R}}

\newcommand{\oo}{\infty}

\newcommand{\x}{\times}

\newcommand{\nin}{\not\in}
\newcommand{\cat}{\hat{\ }}
\newcommand{\sbst}{\subseteq}
\newcommand{\spst}{\supseteq}
\newcommand{\sm}{\setminus}

\newcommand{\op}{\operatorname}

\newcommand{\cov}{\op{cov}}

\newcommand{\non}{\op{non}}

\long\def\note#1\endnote%
{\ifNote{\par\medskip\tt Note: #1}\medskip\par
                         \else{}
                         \fi}

\newtheorem{thm}{Theorem}[section]
\newtheorem{cor}[thm]{Corollary}

\newtheorem{lem}[thm]{Lemma}
\newtheorem{prop}[thm]{Proposition}
\newtheorem{prob}[thm]{Problem}

\theoremstyle{remark}

\newtheorem{rem}[thm]{Remark}

\theoremstyle{definition}
\newtheorem{defn}[thm]{Definition}

\newcommand{\be}{\begin{enumerate}}
\newcommand{\ee}{\end{enumerate}}
\newcommand{\bi}{\begin{itemize}}
\newcommand{\ei}{\end{itemize}}
\newcommand{\itm}{\item}

\title[The Baer-Specker group]{The combinatorics of the Baer-Specker group}
\author{Micha\l{} Machura}
\address{Institute of Mathematics, University of Silesia, ul.\ Bankowa 14, 40-007 Katowice, Poland}
\email{machura@ux2.math.us.edu.pl}

\author{Boaz Tsaban}
\address{Department of Mathematics,
Weizmann Institute of Science,
Rehovot 76100, Israel;
and
Department of Mathematics,
Bar-Ilan University,
Ramat-Gan 52900, Israel}
\email{tsaban@math.biu.ac.il}
\urladdr{http://www.cs.biu.ac.il/\~{}tsaban}
\thanks{The first author was supported by the EU Research and Training Network HPRN-CT-2002-00287.
The second author was supported by the Koshland Center for Basic Research.}

\begin{document}

\begin{abstract}
We study subgroups of $\Bgp$ which possess group theoretic properties
analogous to properties introduced by Menger (1924), Hurewicz (1925),
Rothberger (1938), and Scheepers (1996).
The studied properties were introduced independently by Ko\v{c}inac and Okunev.
We obtain purely combinatorial characterizations of these properties, and combine
them with other techniques to solve several questions of Babinkostova, Ko\v{c}inac, and Scheepers.
\end{abstract}

\maketitle

\section{Introduction}

The groups $\Z^k$ ($k\in\N$) are discrete, and the
classification up to isomorphism of their (topological) subgroups
is trivial. But already for the countably infinite power $\Bgp$ of $\Z$, the situation
is different. Here the product topology is nontrivial,
and the subgroups of $\Bgp$ make a rich source of examples
of non-isomorphic topological groups (see the papers in \cite{BlAb}, and
references therein).
Because of the early works of Baer \cite{Baer37} and Specker \cite{Spec50},
the group $\Bgp$ is commonly named \emph{the Baer-Specker group}.

We study properties of subgroups of the Baer-Specker group,
which are preserved under continuous group homomorphisms.
While the definitions of these properties contain a topological
ingredient, they all turn out to be equivalent to purely
combinatorial properties.
One of these properties is of special interest in light of
classical conjectures of Menger and Hurewicz.

This paper is organized as follows. In Section
\ref{Gdefs} we define the studied group theoretic properties.
These definitions, due to Okunev and Ko\v{c}inac independently,
involve coverings of the group by translates of open sets.
In Section \ref{comb} we provide combinatorial characterizations for
these properties.
In Section \ref{revisited}, we discuss the Hurewicz Conjecture for groups,
in light of a recent result of Babinkostova.
In Section \ref{contimages} we generalize the properties for arbitrary sets,
prove their preservation under uniformly continuous images, and give a characterization
of each of the general topological properties in terms of the group theoretic ones.
In Section \ref{critcards} we describe the minimal cardinalities of counter-examples
for the group theoretic properties.
Section \ref{constructions1} shows that the group theoretic properties do not coincide,
and Section \ref{constructions2} shows that none of the group theoretic properties
coincides with its general topological counterpart,
answering several questions of Babinkostova, Ko\v{c}inac, and Scheepers posed in \cite{coc11}.
In Section \ref{constructions3} we give a systematic list of constructions witnessing the diversity
of subgroups of the Baer-Specker group.
An informal thesis that emerges from that section is that $\Bgp$ is ``universal'' for boundedness properties
of groups.

\section{Boundedness notions for groups}\label{Gdefs}

\begin{defn}\label{bddgps}
Assume that $(G,\cdot)$ is a topological group.
\be

\itm $G$ is \emph{Menger-bounded} if for each sequence
$\seq{U_n}$ of neighborhoods of the unit, there exist finite sets
$F_n\sbst G$, $n\in\N$, such that $G=\Union_n F_n\cdot U_n$.\footnote{Throughout
the paper, $A\cdot B$ stands for $\{a\cdot b : a\in A,\ b\in B\}$, and
$a\cdot B$ stands for $\{a\cdot b : b\in B\}$.}

\itm $G$ is \emph{Scheepers-bounded} if for each sequence
$\seq{U_n}$ of neighborhoods of the unit, there exist finite sets
$F_n\sbst G$, $n\in\N$, such that for each finite set $F\sbst G$,
there is $n$ such that $F\sbst F_n\cdot U_n$.

\itm $G$ is \emph{Hurewicz-bounded} if for each sequence
$\seq{U_n}$ of neighborhoods of the unit, there exist finite sets
$F_n\sbst G$, $n\in\N$, such that for each $g\in G$, $g\in F_n\cdot U_n$
for all but finitely many $n$.

\itm $G$ is \emph{Rothberger-bounded} if for each sequence
$\seq{U_n}$ of neighborhoods of the unit, there exist elements
$a_n\in G$, $n\in\N$, such that $G=\Union_n a_n\cdot U_n$.

\ee
\end{defn}

Several instances of these properties were studied in, e.g.,
\cite{TkaIntro, Hernandez, HRT, KMexample, BNS, o-bdd}.
A study from a more general point of view was initiated in \cite{Koc03, coc11, Bab05}.
These properties are obtained from the following
general topological properties by
restricting attention to open covers of the form $\{a\cdot U : a\in G\}$,
where $U$ is an open neighborhood of the unit.

\begin{defn}\label{sps}
Assume that $X$ is a topological space.
\be

\itm $X$ has the \emph{Menger property} \cite{Menger24}
if for each sequence
$\seq{\cU_n}$ of open covers of $X$, there exist finite sets
$\cF_n\sbst\cU_n$, $n\in\N$, such that $\Union_{n\in\N} \cF_n$ is a cover of $X$.

\itm $X$ has the \emph{Scheepers property} \cite{coc1}
if for each sequence
$\seq{\cU_n}$ of open covers of $X$, there exist finite sets
$\cF_n\sbst\cU_n$, $n\in\N$, such that
for each finite set $F\sbst X$,
there is $n$ such that $F\sbst \Union_{U\in\cF_n} U$.

\itm $X$ has the \emph{Hurewicz property} \cite{Hure25, Hure27}
if for each sequence
$\seq{\cU_n}$ of open covers of $X$, there exist finite set
$\cF_n\sbst\cU_n$, $n\in\N$, such that
for each element $x\in X$,
$x\in\Union_{U\in\cF_n} U$ for all but finitely many $n$.

\itm $X$ has the \emph{Rothberger property} \cite{Roth38}
if for each sequence
$\seq{\cU_n}$ of open covers of $X$, there exist elements
$U_n\in\cU_n$, $n\in\N$, such that
$X=\Union_{n\in\N} U_n$.
\ee
\end{defn}
Except for the second, all these properties are classical.
They share the same structure and can be defined in a unified manner \cite{coc1, coc2}.
These properties were analyzed in many papers
and form an active area of mathematical research -- see \cite{LecceSurvey, KocSurv, ict, BZSPM}
and references therein.

The group theoretic properties are the main object of study in this paper,
but we also consider their relation to the general topological properties.
Clearly, the group theoretic properties are related as follows:
{\small
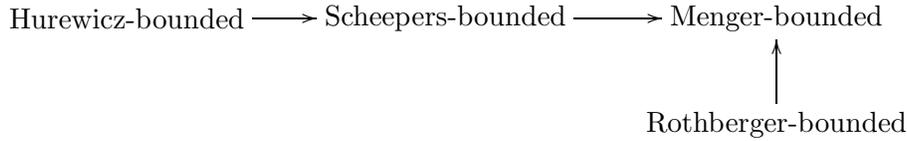
\begin{center}
$\xymatrix{
\mbx{Hurewicz-bounded} \ar[r]        & \mbx{Scheepers-bounded} \ar[r] & \mbx{Menger-bounded}\\
 & & \mbx{Rothberger-bounded}\ar[u]
}$
\end{center}
}
In addition, they are all hereditary for subgroups and preserved under continuous homomorphisms
\cite{Hernandez, coc11, Bab05}.
All properties in the top row hold for $\sigma$-compact groups, and
therefore for subgroups of $\sigma$-compact groups.
In particular, all properties in the top row hold for subgroups of the \emph{Cantor group} $\Cgp$.
As we shall see, the situation is quite different in the case of the Baer-Specker group $\Bgp$.

\section{Purely combinatorial characterizations and some consequences}\label{comb}

We use the convention that $0\in\N$.
For the sake of clarity, we use the following self-evident notations.
$\NNup$ is the collection of all strictly increasing elements of $\NN$, and
$\FinSeqs{\Z}$ is the collection of all finite sequences of integers.
The canonical basis for the topology of $\Bgp$ consists of the sets
$$\bo{s} = \{f\in\Bgp : s\sbst f\}$$
where $s$ ranges over $\FinSeqs{\Z}$.
For natural numbers $k<m$, $[k,m)=\{k,k+1,\dots,m-1\}$,
and $[k,\oo)=\{k,k+1,\dots\}$.
For a partial function $f:\N\to\Z$, $|f|$ is the function with the same domain,
which satisfies $|f|(n)=|f(n)|$, where in this case $|\cdot|$ denotes the absolute value.
For partial functions $f,g:\N\to\Z$, $f\le g$ means: For each $n$ in the domain of $f$,
$n$ is also in the domain of $g$, and $f(n)\le g(n)$.
Similarly, $f\le k$ means: For each $n$ in the domain of $f$,
$f(n)\le k$.
The quantifiers $(\exists^\oo n)$ and $(\forall^\oo n)$ stand for ``there exist infinitely many $n$'' and
``for all but finitely many $n$'', respectively.
Finally, the identity element of $\Bgp$, the constantly zero sequence, is also denoted by $0$.

\subsection{Menger-bounded groups}

\begin{thm}\label{MenBdd}
Assume that $G$ is a subgroup of $\Bgp$. The following
conditions are equivalent:
\be
\itm $G$ is Menger-bounded.
\itm For each $h\in\NNup$, there is $f\in\NN$
such that:
$$(\forall g\in G)(\exists n)\ |g|\rest [0,h(n))\le f(n).$$
\itm For each $h\in\NNup$, there is $f\in\NN$
such that:
$$(\forall g\in G)(\exists^\oo n)\ |g|\rest [0,h(n))\le f(n).$$
\itm There is $f\in\NN$ such that:
$$(\forall g\in G)(\exists^\oo n)\ |g|\rest [0,n)\le f(n).$$
\ee
\end{thm}
\begin{proof}
$(1\Impl 2)$ Fix $h \in \NNup$.
For each $n$, take $U_n = \bo{0\rest [0,h(n))}$ and find finite $F_n\sbst G$ such that
$G\sbst\Union_n (F_n + U_n)$.
Define $f\in\NN$ by
$$f(n) = \max\{|a(k)| : a\in F_n\mbox{ and }k<h(n)\}$$
for each $n$.
For each $g\in G$ there are $n$ and $a\in F_n$ such that $g\in a+U_n = \bo{a\rest[0,h(n))}$,
that is, $g\rest[0,h(n))=a\rest[0,h(n))$. Thus, $|g|\rest[0,h(n))=|a|\rest[0,h(n))\le f(n)$
for this $n$.

$(2\Impl 1)$
Assume that $\seq{U_n}$ is a sequence of neighborhoods of $0$ in $\Bgp$.
Take $h\in\NNup$ such that $\bo{0\rest[0,h(n))}\sbst U_n$ for each $n$.
Apply (2) for $h$ to obtain $f$.
For each $n$ and each $s\in\Z^{[0,h(n))}$ with $|s|\le f$,
choose (if possible) $a_s \in G$ such that
$a_s\rest[0,h(n)) = s$. If this is impossible, take $a_s = 0$.
Let $F_n = \{a_s : s\in\Z^{[0,h(n))},\ |s|\le f\}$.
We claim that $G\sbst\Union_n(F_n+U_n)$.
For each $g\in G$, there is $n$ such that $|g|\rest[0,h(n)) \le f(n)$,
and therefore there is $s\in\Z^{[0,h(n))}$ such that $g\rest[0,h(n))=s=a_s\rest[0,h(n))$,
thus
$$g\in \bo{a_s\rest[0,h(n))} = a_s+\bo{0\rest[0,h(n))}\sbst a_s+U_n\sbst F_n+U_n.$$

$(2\Impl 3)$ This is achieved by partitioning $\N$ to infinitely many
infinite pieces and applying the arguments in $(1\Impl 2)$ to each piece separately.

$(3\Impl 2)$ and $(3\Impl 4)$ are trivial.

$(4\Impl 3)$ This was pointed out by Banakh and Zdomskyy, and later independently
by Simon. Indeed, fix any $h\in\NNup$. Let $f$ be as in $(4)$.
We may assume that $f$ is increasing. Define $\tilde f(n)=f(h(n+1))$ for each $n$.
Fix $g\in G$.
For each $n>h(0)$ with $|g|\rest [0,n)\le f(n)$, let $m$ be such that $n\in\intvl{h}{m}$.
Then $|g|\rest [0,h(m))\le |g|\rest [0,n) \le f(n) \le f(h(m+1)) = \tilde f(m)$.
There are infinitely many such $n$'s, and therefore infinitely many such $m$'s.
\end{proof}

\subsection{Scheepers-bounded groups}

\begin{thm}\label{SchBdd}
Assume that $G$ is a subgroup of $\Bgp$. The following
conditions are equivalent:
\be
\itm $G$ is Scheepers-bounded.
\itm For each $h\in\NNup$, there is $f\in\NN$
such that:
$$(\forall\mbox{finite }F\sbst G)(\exists n)(\forall g\in F)\ |g|\rest [0,h(n))\le f(n).$$
\itm For each $h\in\NNup$, there is $f\in\NN$
such that:
$$(\forall\mbox{finite }F\sbst G)(\exists n)(\forall g\in F)\ g\rest [0,h(n))\le f(n).$$
\itm There is $f\in\NN$ such that:
$$(\forall\mbox{finite }F\sbst G)(\exists^\oo n)(\forall g\in F)\ |g|\rest [0,n)\le f(n).$$
\itm There is $f\in\NN$ such that:
$$(\forall\mbox{finite }F\sbst G)(\exists^\oo n)(\forall g\in F)\ g\rest [0,n)\le f(n).$$
\ee
Moreover, in (2) and (3) the quantifier $(\exists n)$ can be replaced
by $(\exists^\oo n)$.
\end{thm}
\begin{proof}
$(3\Impl 2)$ Given a finite $F\sbst G$, apply (3) to the finite set $F\cup -F = \{\pm a : a\in F\}$.

$(5\Impl 4)$ is identical. The remaining implications are proved as in Theorem \ref{MenBdd}.
\end{proof}

\subsection{Hurewicz-bounded groups}
\begin{defn}
A partial ordering $\le^*$ is defined on $\NN$ by:
$f\le^* g$ if $f(n)\le g(n)$ for all but finitely many $n$.
A subset $X$ of $\Bgp$ is \emph{$\le^*$-bounded} if there is $f\in\NN$ such that
for each $g\in X$, $|g|\le^* f$.
\end{defn}

\begin{thm}\label{HurBdd}
Assume that $G$ is a subgroup of $\Bgp$. The following
conditions are equivalent:
\be
\itm $G$ is Hurewicz-bounded.
\itm For each $h\in\NNup$, there is $f\in\NN$
such that:
$$(\forall g\in G)(\forall^\oo n)\ |g|\rest [0,h(n))\le f(n).$$
\itm $G$ is $\le^*$-bounded.
\ee
\end{thm}
\begin{proof}
$(1\Iff 2)$ this is similar to the proof of Theorem \ref{MenBdd}.

$(2\Impl 3)$ is trivial.

$(3\Impl 2)$ Assume that $f_0\in\NN$ witnesses that $G$ is $\le^*$-bounded.
We may assume that $f_0$ is increasing.
For each $g\in G$ let $n_0$ be such that
$|g|\rest [n_0,\oo)\le f_0$, and let $m=\max|g|\rest[0,n_0)$.
Choose $n_1>n_0$ such that $m\le f_0(n_1)$.
Then for each $n\ge n_1$, we have the following:
\begin{eqnarray*}
|g|\rest[0,n_0) \le m\le f_0(n_1)\le f_0(n)\\
|g|\rest[n_0,n) \le f_0\rest[n_0,n)\le f_0(n)
\end{eqnarray*}
In other words, $|g|\rest[0,n)\le f_0(n)$ for each $n\ge n_1$.
Given $h\in\NNup$, define $f\in\NN$ by $f(n)=f_0(h(n))$ for each $n$.
For each large enough $n$, $|g|\rest[0,h(n))\le f_0(h(n))=f(n)$.
\end{proof}

Item (3) in Theorem \ref{HurBdd} is familiar to algebraists under the name
\emph{bounded growth type} (see \cite{Spec50, GoWa80}), and to topologists
as ``subsets of a $\sigma$-compact set.''

\subsection{Rothberger-bounded groups and strong measure zero}
According to Borel \cite{Borel},
a metric space $(X,d)$ has \emph{strong measure zero} if for each sequence $\seq{\epsilon_n}$
of positive reals, there exists a cover $\setseq{U_n}$ of $X$ such that for each $n$,
the diameter of $U_n$ is smaller than $\epsilon_n$. The topology on the Baer-Specker group $\Bgp$
is induced by the metric
$$d(x,y) =
\begin{cases}
\frac{1}{N(x,y)+1} & x \neq y\\
0 & x=y
\end{cases}$$
where for distinct $x,y\in\Bgp$, $N(x,y)=\min\{n : x(n)\neq y(n)\}$.

The equivalence $(2\Iff 3)$ in the following theorem (with $\Cgp$ instead of $\Bgp$)
is from Bartoszy\'nski-Judah \cite{barju},
and the equivalence $(1\Iff 3)$ follows from a general result of
Babinkostova, Ko\v{c}inac, and Scheepers \cite{coc11}.

\begin{thm}[\cite{barju, coc11}]\label{RothBdd}
Assume that $G$ is a subgroup of $\Bgp$. The following
conditions are equivalent:
\be
\itm $G$ is Rothberger-bounded.
\itm For each $h\in\NNup$, there is $\vphi:\N\to\FinSeqs{\Z}$
such that:
$$(\forall g\in G)(\exists n)\ g\rest [0,h(n)) = \vphi(n).$$
\itm $G$ has strong measure zero.
\ee
Moreover, in (2) the quantifier $(\exists n)$ can be replaced
by $(\exists^\oo n)$.
\end{thm}
\begin{proof}
$(1\Impl 2)$ Fix $h \in \NNup$.
For each $n$, take $U_n = \bo{0\rest [0,h(n))}$ and find $a_n \in G$ such that
$G\sbst\Union_n (a_n + U_n)$.
For each $g\in G$ there is $n$ such that $g\in a_n +U_n = \bo{a_n\rest[0,h(n))}$,
that is, $g\rest[0,h(n))=a_n\rest[0,h(n))$. Take $\vphi(n)=a_n\rest[0,h(n))$ for each $n$.

$(2\Impl 1)$
Assume that $\seq{U_n}$ is a sequence of neighborhoods of $0$ in $\Bgp$.
Take $h\in\NNup$ such that $\bo{0\rest[0,h(n))}\sbst U_n$ for each $n$.
Apply (2) for $h$ to obtain $\vphi$.
For each $n$ choose (if possible) $a_n \in G$ such that
$a_n\rest[0,h(n)) = \vphi(n)$. In this is impossible, take $a_n = 0$.
We claim that $G\sbst\Union_n(a_n+U_n)$.
Indeed, for each $g\in G$, there is $n$ such that $g\rest[0,h(n)) = \vphi(n) = a_n\rest[0,h(n))$,
and therefore
$$g\in\bo{a_n\rest[0,h(n))} = a_n+\bo{0\rest[0,h(n))}\sbst a_n+U_n.$$

$(1\Impl 3)$ Given $\seq{\epsilon_n}$, choose for each $n$ a neighborhood of the identity whose
diameter is smaller than $\epsilon_n$. Apply (1) and the fact that the metric on $\Bgp$ is
translation invariant.

$(3\Impl 2)$
Let $h\in\NNup$. For each $n$, take $\epsilon_n=1/h(n)$. By (3), there is a cover $\setseq{U_n}$ of $G$
such that for each $n$, the diameter of $U_n$ is smaller than $\epsilon_n$.
Consequently, each $U_n$ is contained in some $\bo{s_n}$ where $s_n\in\Z^{h(n)}$.
Take $\vphi(n)=s_n$ for each $n$.

The last assertion can be proved as in Theorem \ref{MenBdd}.
\end{proof}

\section{The Hurewicz and Menger Conjectures revisited}\label{revisited}

The notion of Menger-bounded groups was introduced by Okunev (under the name \emph{$o$-bounded groups})
with the aim to have an inner characterization of subgroups of $\sigma$-compact groups.
In the general topological case, this approach goes back to Menger \cite{Menger24},
who conjectured that for metric spaces,
the Menger property (Definition \ref{sps}(1)) characterizes $\sigma$-compactness.
Following that, Hurewicz \cite{Hure25} made the weaker conjecture that the Hurewicz property
(Definition \ref{sps}(3)) characterizes $\sigma$-compactness.
These conjectures turn out to be false \cite{FM, coc2} (generalized in \cite{SFH}).
However, the conjectures also make sense in the group theoretic case.
Since the group theoretic properties are hereditary for subgroups,
we should restate them in the following manner.

\begin{defn}\label{MHC}
The \emph{Menger Conjecture for groups} (respectively, \emph{Hure\-wicz Conjecture for groups})
is the assertion that each metrizable group $G$,
which is Menger-bounded (respectively, Hurewicz-bounded),
is a subgroup of some $\sigma$-compact group.
\end{defn}

Theorem \ref{HurBdd} implies that
the Hurewicz Conjecture is true when restricted to subgroups of $\Bgp$.
A result of Babinkostova implies that it is true for general metrizable groups.
Call a subset $B$ of a topological group $G$ \emph{left totally bounded}
(respectively, \emph{right totally bounded}, \emph{totally bounded})
in $G$ if for each neighborhood $U$ of the identity, there is a finite $F\sbst G$ such
that $B \sbst F\cdot U$ (respectively, $B \sbst U\cdot F$, $B \sbst F\cdot U\cap U\cdot F$).

\begin{thm}[Babinkostova \cite{Bab05}]\label{tb}
For metrizable groups $G$: $G$ is Hurewicz-bounded if, and only if,
$G$ is a union of countably many left totally bounded sets.
\end{thm}

For completeness, we give a direct proof of this result.

\begin{proof}
The ``if'' part is easy. We prove the ``only if'' part.
Fix a metric on $G$ and for each $n$,
let $U_n$ be the ball of radius $1/(n+1)$ centered at the
identity.
Choose finite sets $F_n\sbst G$, $n\in\N$, such that for each $g\in G$
and all but finitely many $n$, $g\in F_n\cdot U_n$.
Then
$$G\sbst \Union_{m\in\N}\bigcap_{n\ge m} F_n\cdot U_n.$$
For each $m$, $\bigcap_{n\ge m} F_n\cdot U_n$ is left totally bounded in $G$.
\end{proof}

A topological group $G$ can be embedded in a $\sigma$-compact group if, and only if,
$G$ is a union of countably many totally bounded sets \cite{TkaIntro}.
The nontrivial implication
follows from the classical fact that the completion of $G$ with respect to its
left-right uniformity is a topological group, and the basic fact that totally bounded and complete
sets are compact.

\begin{cor}\label{HC}
The Hurewicz Conjecture for groups is true.
\end{cor}
\begin{proof}
Assume that $G$ is Hurewicz-bounded. By Theorem \ref{tb},
$G=\Union_n B_n$ with each $B_n$ left totally bounded.
Then
$G = G\inv = \Union_n B_n\inv$,
where $A\inv$ means $\{a\inv : a\in A\}$.
If $B$ is left totally bounded, then $B\inv$ is right totally bounded.
(For each neighborhood $U$ of the identity, take a neighborhood $V$ of the identity
such that $V\inv\sbst U$, and finite $F$ such that $B\sbst F\cdot V\inv$.
Then $B\inv\sbst V\cdot F\inv$.)
Thus,
$$G = \Union_{m,n\in\N} B_n\cap B_m\inv,$$
and each $B_n\cap B_m\inv$ is totally bounded.
That is, $G$ is a union of countably many totally bounded sets.
\end{proof}

On the other hand, by the forthcoming Theorem \ref{CPthm},
already for subgroups of $\Bgp$ the Menger Conjecture is false.

\begin{rem}
Partial results in the direction of Theorem \ref{tb} and Corollary \ref{HC},
and results similar to these, were previously obtained in
\cite{Ban02, MichPhD, BNS} (this list may be nonexhaustive).
\end{rem}

\section{Continuous images}\label{contimages}

The group theoretic properties are preserved under continuous homomorphisms, while
their topological counterparts are preserved under arbitrary continuous functions.
The combinatorial properties characterizing the group theoretic properties
make sense for arbitrary subsets of $\Bgp$ (or of the \emph{Baire space} $\NN$),
and it turns out that they are preserved under uniformly continuous images in $\Bgp$.
A consequence of Specker's work is that every endomorphism of
$\Bgp$ is continuous, and therefore \emph{uniformly} continuous.

\begin{defn}\label{SetBdd}
Abusing terminology, we say that a subset $X$ of $\Bgp$ is Menger-
(respectively, Scheepers-, Hurewicz-, Rothberger-) bounded
if it satisfies the property in Theorem \ref{MenBdd}(4) (respectively, \ref{SchBdd}(5), \ref{HurBdd}(3),
\ref{RothBdd}(2)).
\end{defn}
As far as uniformly continuous images are concerned, we can equivalently work in $\NN$
(the natural homeomorphism from $\Bgp$ to $\NN$ is uniformly continuous in both directions).
In this case, the absolute values in the definitions are not needed.

We say that $Y$ is a \emph{uniformly continuous image} of $X$ if
it is the image of a uniformly continuous function $\Psi:X\to\Bgp$.

\begin{lem}\label{unifcont}
Each of the properties in Definition \ref{SetBdd} is preserved under uniformly continuous images.
\end{lem}
\begin{proof}
To give a unified proof, we use the equivalent formulations of the properties which
use the additional quantifier ``there is $h\in\NNup$''.
Assume that $\Psi:X\to Y$ is a uniformly continuous surjection.
Let $h\in\NNup$ be given for $Y$. As $\Psi$ is uniformly continuous, there exists
for each $n$ a number $h'(n)$ such that for each $x\in X$, $x\rest[0,h'(n))$ determines
$\Psi(x)\rest [0,h(n))$.

Assume that $X$ is Menger-bounded (the remaining cases are similar).
Then there is $f'$ for $h'$ as in \ref{MenBdd}(2).
For each $n$, let $S_n$ be the set of all $s\in \Z^{h'(n)}$ such that $|s|\le f'(n)$,
and such that there is $x_s\in X$ such that $x_s\rest [0,h'(n))=s$.
For each $s\in S_n$, let $r_s = \Psi(x_s)\rest[0,h(n))$ (note that $r_s$ depends only on $s$).
Define $f\in\NN$ by
$$f(n)=\max\{|r_s(k)| : s\in S_n,\ k<h'(n)\}$$
for each $n$ (if $S_n=\emptyset$ take $f(n)=0$).
Then $f$ is as required in \ref{MenBdd}(2) for $Y$ and $h$.
\end{proof}

\begin{thm}\label{morechars}
Assume that $X\sbst\Bgp$.
$X$ is Menger- (respectively, Scheepers-, Hurewicz-, Rothberger-) bounded, if, and only if,
for each uniformly continuous image $Y$ of $X$ in $\Bgp$ and each $h\in\NNup$, there is $f\in\NN$ such that:
\be
\itm $(\forall y\in Y)(\exists n)\ |y|\rest \intvl{h}{n}\le f\rest\intvl{h}{n}$,
\itm $(\forall\mbox{finite }F\sbst Y)(\exists n)(\forall y\in F)\ |y|\rest \intvl{h}{n}\le f\rest\intvl{h}{n}$,
\itm $(\forall y\in Y)(\forall^\oo n)\ |y|\rest \intvl{h}{n}\le f\rest\intvl{h}{n}$,
\itm $(\forall y\in Y)(\exists n)\ y\rest \intvl{h}{n} = f\rest\intvl{h}{n}$,
\ee
respectively.
Moreover, in each of the above items the quantifier $(\exists n)$ can be replaced
by $(\exists^\oo n)$.
\end{thm}
\begin{proof}
By Lemma \ref{unifcont}, the implication from the definitions to the new characterizations is immediate.
We prove the converse direction.

(1) Given $h\in\NNup$ for $X$, define $h'\in\NNup$ by
$h'(0)=0$, and $h'(n)=h(0)+\dots+h(n-1)$ for each $n>0$.
Define $\Psi:X\to\Bgp$ by
$$\Psi(x)(k) = x(k-h'(n))$$
whenever $k\in\intvl{h'}{n}$.
$\Psi$ is uniformly continuous (in fact, it is a homomorphism).
Let $Y$ be the image of $\Psi$, and take $f'$ as in (1) for $h'$.
Then for each $x\in X$, there is $n$ such that
$x(k-h'(n))=\Psi(x)(k)\le f'(k)$ whenever $k\in\intvl{h'}{n}$.
Define $f(n)=\max\{f'(k) : k\in\intvl{h'}{n}\}$ for each $n$.
Then $f$ is as in \ref{MenBdd}(2) for $X$ and $h$.

(2) is similar, and (3) is trivial.

(4) Argue as in (1), and define $\vphi(n)(k)= f'(k+h'(n))$ for each $k<h(n)$.
\end{proof}

\begin{rem}
The proof of Theorem \ref{morechars} actually shows that in the context of groups,
the same assertions hold when restricting attention to endomorphisms rather than arbitrary
uniformly continuous functions.
\end{rem}

We now turn to arbitrary continuous images.

\begin{thm}[\cite{Hure25, Rec94, huremen1}]\label{rec}
Assume that $X\sbst\Bgp$.
$X$ has the Menger (respectively, Scheepers, Hurewicz, Rothberger) property if, and only if,
for each continuous image $Y$ of $X$ in $\Bgp$, there is $f\in\NN$ such that:
\be
\itm $(\forall y\in Y)(\exists n)\ |y|(n)\le f(n)$,
\itm $(\forall\mbox{finite }F\sbst Y)(\exists n)(\forall y\in F)\ |y|(n)\le f(n)$,
\itm $(\forall y\in Y)(\forall^\oo n)\ |y|(n)\le f(n)$,
\itm $(\forall y\in Y)(\exists n)\ y(n) = f(n)$,
\ee
respectively.
Moreover, in each of the above items the quantifier $(\exists n)$ can be replaced
by $(\exists^\oo n)$.
\end{thm}

In light of Theorem \ref{RothBdd}, we have that the following corollary
extends Fremlin-Miller's Theorem 1 from \cite{FM} (cf.\ \cite{BZSPM}).

\begin{thm}\label{FMlike}
Assume that $X\sbst\Bgp$.
$X$ has the Menger (respectively, Scheepers, Hurewicz, Rothberger) property if, and only if,
each continuous image $Y$ of $X$ in $\Bgp$
is Menger- (respectively, Scheepers-, Hurewicz-, Rothberger-) bounded.
\end{thm}
\begin{proof}
$(\Impl)$ Each of the topological properties is preserved under continuous images and implies the
corresponding group theoretic property.

$(\Leftarrow)$ We treat the Menger case, the other cases being similar.
If $Y$ is a continuous image of $X$ in $\Bgp$,
then by the assumption $Y$ is Menger-bounded. By Theorem \ref{MenBdd}(2) for $Y$,
we have, in particular, that (1) of Theorem \ref{rec} holds for $Y$.
\end{proof}

These results and those in the coming sections imply that none of the group-theoretic
properties (considered for general subsets of $\Bgp$) is preserved under continuous images.

\section{Critical cardinalities}\label{critcards}

A subset $D$ of $\NN$ is \emph{dominating} if
for each $g\in\NN$ there exists $f\in D$ such that $g\le^* f$.
Let $\fb$ denote the minimal cardinality of a $\le^*$-unbounded
subset of $\NN$, and $\fd$ denote the minimal cardinality of a dominating
subset of $\NN$. In addition, let $\cov(\cM)$ denote the minimal cardinality
of a subset $Y$ of $\NN$ such that there is no $f\in\NN$ such that
for each $y\in Y$, $f(n)=y(n)$ for some $n$.
The name $\cov(\cM)$ comes from the fact that this cardinal is also the minimal
cardinality of a cover of $\R$ by meager sets \cite{barju}.
Finally, let $\fc = 2^{\alephes}$ denote the cardinality of the continuum.

The mentioned cardinals are related as follows, where an arrow means $\le$:
\begin{center}
$\xymatrix{
\fb \ar[r] & \max\{\fb,\cov(\cM)\}\ar[r] & \fd \ar[r] & \fc\\
\aleph_1 \ar[u]\ar[r] & \cov(\cM)\ar[u]
}$
\end{center}
No additional (weak) inequalities among these cardinals can be proved, see \cite{BlassHBK}.

The \emph{critical cardinality} of a property $P$ of subsets of $\Bgp$ is:
$$\non(P)=\min\{|X| : X\sbst\Bgp\mbox{ and }X\mbox{ does not satisfy }P\}.$$
Consider the properties (1)--(4) in Theorem \ref{rec}.
The critical cardinality of (1) is $\fd$.
It is not difficult to see that $\fd$ is also the critical cardinality of (2).
The critical cardinality of (3) is $\fb$,
and $\cov(\cM)$ is the critical cardinality of (4).

\begin{cor}\label{crits}
\be
\itm The critical cardinalities of the properties Menger, Menger-bounded (for sets or groups),
Scheepers, and Scheepers-bounded (for sets or groups) are all $\fd$.
\itm The critical cardinalities of the properties Hurewicz and Hurewicz-bounded (for sets or groups) are
$\fb$.
\itm The critical cardinalities of the properties Rothberger and Rothberger-bounded (for sets or groups) are
$\cov(\cM)$.
\ee
\end{cor}
\begin{proof}
For the topological properties this follows from Theorem \ref{rec}.
Consequently, Theorem \ref{FMlike} implies the assertions for the bounded version of the properties for sets.
To get the property for groups, take a witness $Y\sbst\Bgp$ for the critical cardinality of the
same property for arbitrary sets, and consider $\langle Y\rangle$ (which has the same cardinality as $Y$).
\end{proof}

\section{Comparison of the group theoretic properties}\label{constructions1}

Could any implication---which is not obtained by composition of existing ones---be
added to the diagram in Section \ref{Gdefs}? We give examples ruling out
almost all possibilities.

The most simple example is the following.

\begin{thm}\label{cantorgroup}
Let $X=\Cantor\sbst\Bgp$. Then $G=\langle X\rangle$ is $\sigma$-compact
(in particular, it is Hurewicz-bounded) but is not Rothberger-bounded.
\end{thm}
\begin{proof}
$X$ is compact and therefore so are all of its finite powers.

\begin{lem}\label{powers}
Let $M$ be a topological group.
Assume that $P$ is a property which is preserved under uniformly continuous images and countable unions.
If all finite powers of a subset $X$ of $M$ have the property $P$, then all finite powers
of the group $G=\langle X\rangle$ have the property $P$.
\end{lem}
\begin{proof}
For each $n$,
$$G_n=\{m_1x_1+\dots + m_nx_n : m_1,\dots,m_n\in\Z,\ x_1,\dots,x_n\in X\}$$
is a union of countably many continuous images of $X^n$,
Thus, for each $k$, $(G_n)^k$ is a union of countably many continuous images of $X^{nk}$,
and therefore has the property $P$.
Thus, for each $k$, $G^k=\langle X\rangle^k = \Union_n(G_n)^k$ has the property $P$.
\end{proof}

As $\sigma$-compactness is preserved under continuous images and countable unions,
we have by Lemma \ref{powers} that $G=\langle X\rangle$ is $\sigma$-compact.
On the other hand, it is clear by Theorem \ref{RothBdd}(2) that $X$ (and in particular $G$)
is not Rothberger-bounded.
\end{proof}

\begin{thm}\label{CPthm}
There exists a Scheepers-bounded subgroup of $\Bgp$ which is not Hurewicz-bounded.
\end{thm}
\begin{proof}
By a theorem of Chaber and Pol \cite{ChaPol} (see also \cite{SFH}),
there is a subset $X$ of $\Bgp$
such that all finite powers of $X$ have the Menger property, but $X$ is not contained
in any $\sigma$-compact subset of $\Bgp$.
Let $G=\langle X\rangle$.

Menger's property is preserved under continuous images and countable unions.
By Lemma \ref{powers}, $G$ has the Menger property.
By the following lemma, $G$ is Scheepers-bounded.

\begin{lem}[\cite{coc11, BZSPM}]\label{Mengerinpowers}
A topological group $G$ is Scheepers-bounded if, and only if,
all finite powers of $G$ are Menger-bounded.\hfill\qed
\end{lem}

Finally, as $G\spst X$, $G$ is not contained in any $\sigma$-compact subset of $\Bgp$.
By Theorem \ref{HurBdd}, $G$ is not Hurewicz-bounded.
\end{proof}

By a deep result of Laver \cite{Laver}, it is consistent that
all strong measure zero sets of reals are countable (and therefore have all properties
considered in this paper).
Thus, special hypotheses are necessary for constructions of nontrivial Rothberger-bounded groups.
One can replace all of the special hypotheses made in this paper concerning equalities among cardinals
by \CH{} or just \MA{}, which are both stronger than our hypotheses.

\begin{thm}\label{RBNotHB}
Assume that $\cov(\cM)=\fc$. Then there exists a subgroup of $\Bgp$ which
is Rothberger-bounded (and Scheepers-bounded) but not Hurewicz-bounded.
\end{thm}
\begin{proof}
For a cardinal $\kappa$, $L\sbst\Bgp$ is a \emph{$\kappa$-Luzin} set if $|L|\ge\kappa$, and
for each meager set $M\sbst \Bgp$, $|L\cap M|<\kappa$.
Our assumption implies the existence of a $\cov(\cM)$-Luzin set $L\sbst\Bgp$ such that
all finite powers of $L$ have the Rothberger property \cite{coc2}.
Take $G=\langle L\rangle$.
As the Rothberger property is preserved under continuous images and countable unions,
we have by Lemma \ref{powers} that all finite powers of $G$ have the Rothberger (in particular, Menger)
property.
It follows from Lemma \ref{Mengerinpowers} that $G$ is Scheepers-bounded.

\begin{lem}
Assume that $L\sbst\Bgp$ is a $\kappa$-Luzin set. Then $L$ is $\le^*$-unbounded.
\end{lem}
\begin{proof}
For each $f\in\NN$, the set $M=\{g\in\Bgp : |g|\le^* f\}$ is meager in $\Bgp$,
thus $|L\cap M|<\kappa=|L|$, and in particular $L\not\sbst M$.
\end{proof}
As $G\spst L$, $G$ is $\le^*$-unbounded.
By Theorem \ref{MenBdd}, $G$ is not Hurewicz-bounded.
\end{proof}

\begin{rem}
Assuming $\cov(\cM)=\fc$, there exist subgroups of $\R^\N$ which are $\cov(\cM)$-Luzin sets
\cite{o-bdd}. However, in $\Bgp$ this is impossible: The subset $2\Bgp=\{2f : f\in\Bgp\}$ of
$\Bgp$ is nowhere dense, but $G\cap 2\Bgp\spst 2G$ has the same cardinality as $G$.
\end{rem}

In summary, as far as the group theoretic properties are concerned,
no implications can be added among those corresponding to the classical Menger, Hurewicz, and Rothberger properties.
Concerning the newer Scheepers property, only the following (related) problems
remain open (compare this to Problems 1 and 2 of \cite{coc2}).

\begin{prob}\label{probs}
\mbox{}
\be
\itm Is every Menger-bounded subgroup of $\Bgp$ Scheepers-bounded?\/\footnote{This was
recently answered in the negative in: M. Machura, S. Shelah, and B. Tsaban,
\emph{Squares of Menger-bounded groups}
(preprint).}
\itm And if not, is every Rothberger-bounded subgroup of $\Bgp$ Scheep\-ers-bounded?
\ee
\end{prob}

\section{The topological properties are strictly stronger}\label{constructions2}

Ko\v{c}inac in \cite{Koc03} and Babinkostova, Ko\v{c}inac, and Scheepers in their
first version of \cite{coc11},
asked whether each group theoretic property coincides with its topological counterpart.
Note that in Section \ref{constructions1}, all groups exhibiting some group theoretic
properties actually exhibited the corresponding topological property.
We show that in general this is not the case.

In some of the cases, it will be easier to carry out our
constructions in $\Cgp$ rather than in $\Bgp$. The extension to
$\Bgp$ is as follows. Assume that $H$ is a subgroup of $\Cgp$, and
$G=\langle H\rangle$ is the subgroup of $\Bgp$ generated by $H$ as a subset
of $\Bgp$. Then $G$ is a subgroup of $\langle \Cantor\rangle$, which by
Theorem \ref{cantorgroup} is $\sigma$-compact. Thus, $G$ is
Hurewicz-bounded. To see that $G$ violates a property when $H$
does, we will use the following.

\begin{lem}\label{powers2}
Assume that $P$ is a property which is preserved under uniformly continuous images,
$H$ is a subgroup of $\Cgp$, and $G=\langle H\rangle$ is the subgroup of $\Bgp$ generated by
$H$ as a subset of $\Bgp$.
Then:
For each $k$, if $H^k$ does not have the property $P$, then $G^k$ does not have the property $P$.
\end{lem}
\begin{proof}
The mapping $\Psi : \Bgp \to\Cgp$ defined by
$$\Psi(f)(n) = f(n) \bmod 2$$
for each $n$ is a continuous group homomorphism (in particular, it is uniformly continuous),
and therefore
$$\Psi[G]=\Psi[\langle H\rangle]=\langle \Psi[H]\rangle=\langle H\rangle=H.$$
Thus, $H^k$ is a continuous homomorphic
image of $G^k$. Thus, if $G^k$ has the property $P$, then so does
$H^k$.
\end{proof}

Identify $\Cgp$ with $\PN$ by taking characteristic functions.
The group operation on $\PN$ induced by this identification
is $\Delta$, the symmetric difference.
The \emph{Rothberger space}, denoted $\roth$, is the subspace
of $\PN$ consisting of all infinite sets of natural numbers.
Let $\Fin\sbst\PN$ denote the collection of finite sets of natural numbers.
One can think of the following construction and some of the others
in the sequel as being carried out in the \emph{Rothberger group}
$\PN/\Fin$. However, we prefer to argue directly.
$\Fin$ is a subgroup of $\PN$.
Let
$$e:\roth\to\NNup$$
denote the continuous function which assigns to each infinite subset of $\N$ its
increasing enumeration.

\begin{thm}\label{HBNotM}
There exists a subgroup $G$ of a $\sigma$-compact subgroup of $\Bgp$
(thus, $G$ is Hurewicz-bounded) which does not have the Menger property.
\end{thm}
\begin{proof}
If $G$ is a subgroup of $\PN$ with these properties then, by
Lemma \ref{powers2} (here $P$ is the Menger property) and
the discussion preceding it, so is the group it generates in $\Bgp$.
Thus, we can work in $\PN$.

\begin{lem}\label{domi}
There is a family $D\sbst\roth$ such
that the subgroup $G=\langle D\rangle$ of $\PN$ is contained
in $\roth\cup\{\emptyset\}$, and $e[D]\sbst\NN$ is
dominating.
\end{lem}
\begin{proof}
Fix a dominating subset $\{f_\alpha : \alpha < \fd\}$ of $\NNup$.
Define $D=\{d_\alpha : \alpha<\fd\}\sbst\roth$ by induction on $\alpha$:
At step $\alpha$,
let
$$G_\alpha=\langle \{d_\beta : \beta<\alpha\}\cup\Fin\rangle = \langle \{d_\beta : \beta<\alpha\}\rangle\Delta\Fin.$$
Then $|G_\alpha|=\max\{|\alpha|,\alephes\}<\fd$.
There are continuum many $g\in\NNup$ such that $f_\alpha\le^* g$,
so we can choose $g_\alpha\in\NNup$ such that $f_\alpha\le^* g_\alpha\nin e[G_\alpha\cap\roth]$.
Take $d_\alpha = e\inv(g_\alpha)=\im(g_\alpha)$.
\end{proof}

Take $D$ as Lemma \ref{domi}, and
let $G=\langle D\rangle$ be the generated group in $\PN$.
We claim that $G$ does not have the Menger property.
We will use the fact that the Menger property is
stable under removing finitely many points.
The Menger property is hereditary for closed subsets
and preserved under countable unions.
\begin{lem}\label{removept}
Assume that a property $P$ is hereditary for closed subsets
and preserved under countable unions.
If $X\sbst\PN$ has the property $P$,
then for each $x\in X$, $X\sm\{x\}$ has the property $P$.
\end{lem}
\begin{proof}
A well known property of the Cantor space (and therefore also of $\PN$) is
that for each $x\in\PN$, $\PN\sm\{x\}$ is an $F_\sigma$ subset of $\PN$.
Let $x\in X$. Then $X\sm\{x\}=\Union_n C_n$ where each $C_n$ is a closed subset of $X$,
and therefore has the property $P$.
Consequently, $\Union_nC_n=X\sm\{x\}$ has this property, too.
\end{proof}
Assume that $G$ has the Menger property.
Then by Lemma \ref{removept}, so does $G\sm\{\emptyset\}$, which is a subset
of $\roth$. Now, $e[D]\sbst e[G\sm\{\emptyset\}]\sbst\NN$ and $e[D]$ is dominating.
Thus $e[G\sm\{\emptyset\}]$, a continuous image  of $G\sm\{\emptyset\}$, is dominating.
This contradicts Theorem \ref{rec}(1).
\end{proof}

We now treat the Rothberger-bounded groups.
Recall the comments made before Theorem \ref{RBNotHB}.

\begin{thm}\label{RBNotMen}
Assume that $\cov(\cM) = \fb =\fd $.
Then there exists a subgroup $G$ of a $\sigma$-compact subgroup of $\Bgp$
such that all finite powers of $G$ are
Rothberger-bounded (and Hurewicz-bounded), but $G$ does not have the Menger property.
\end{thm}
\begin{proof}
We first prove the assertion for subgroups of $\PN$.

Construct $D\sbst\roth$ as in Lemma \ref{domi},
but this time using the fact that at each stage $\alpha$,
$|\{g_\beta : \beta<\alpha\}|<\fb$ to make sure that for each $\alpha<\beta$,
$g_\beta\le^* g_\alpha$.

Having this, set $X = D \cup \Fin$.
As $\fb\le\cov(\cM)$, we have by Corollary 14 of \cite{ideals} (see \cite{SFH})
that all finite powers of $X$ have the Rothberger property.
By Lemma \ref{powers}, the subgroup $G=\langle X\rangle$ of $\PN$ (and all its finite powers)
has the Rothberger property.
In particular, it is Rothberger-bounded.
Consequently, so is its subgroup $H=\langle D\rangle$.
In the proof of Theorem \ref{HBNotM} it was shown that
$H$ cannot have the Menger property.

In $\Bgp$, take $\tilde G=\langle G\rangle$ and $\tilde H=\langle H\rangle$.
$H,G$ are subgroups of the $\sigma$-compact group $\langle \Cantor\rangle$.
By Lemma \ref{powers}, all finite powers of $\tilde G$ have
the Rothberger property. In particular they are Rothberger-bounded and
thus the same holds for the subgroup $\tilde H$ of $\tilde G$.
By Lemma \ref{powers2}, $\tilde H$ does not have the Menger property.
\end{proof}

\begin{rem}\label{RBNotM2}
The hypothesis $\fb\leq\cov(\cM)$ suffices to obtain a subgroup
$G$ of $\Bgp$ such that all finite powers of $G$ are
Rothberger-bounded (and Hurewicz-bounded), but $G$ does not have the Hurewicz property.
To achieve that, we only make sure in the last proof that $D$ is $\le^*$-unbounded, without
necessarily having it dominating.
\end{rem}

\section{Additional examples}\label{constructions3}

The immediate relationships among the properties considered in this paper
are summarized in the following diagram.

{\tiny
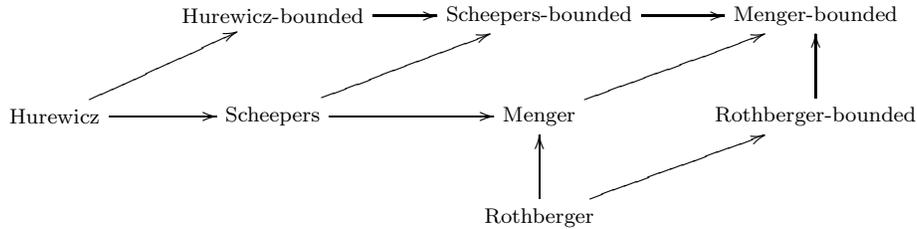
\begin{center}
$\xymatrix{
                      & \mbx{Hurewicz-bounded} \ar[r]        & \mbx{Scheepers-bounded} \ar[r] & \mbx{Menger-bounded}\\
\mbx{Hurewicz} \ar[r]\ar[ur] & \mbx{Scheepers} \ar[r]\ar[ur] & \mbx{Menger}\ar[ur] & \mbx{Rothberger-bounded}\ar[u]\\
                      &                       & \mbx{Rothberger}\ar[u]\ar[ur]
}$
\end{center}
}

An interesting question is, whether additional relationships can be proved to hold.
In other words: Which of the settings consistent with the present implications
can be realized by a single group? As Problem \ref{probs} is still open,
we can only hope to settle the possible settings in the collapsed diagram where
``Scheepers'' and ``Menger'' are identified, that is, consider only the
classical properties:

{\small
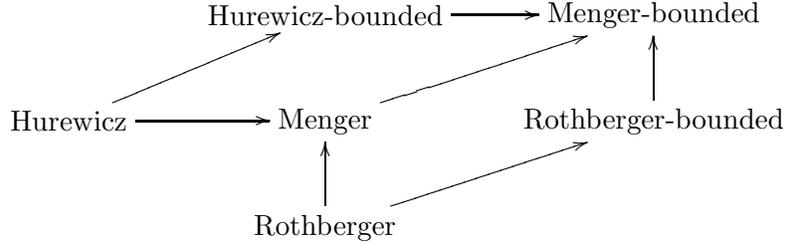
\begin{center}
$\xymatrix{
                      & \mbx{Hurewicz-bounded} \ar[r]        & \mbx{Menger-bounded}\\
\mbx{Hurewicz} \ar[r]\ar[ur] & \mbx{Menger}\ar[ur] & \mbx{Rothberger-bounded}\ar[u]\\
                       & \mbx{Rothberger}\ar[u]\ar[ur]
}$
\end{center}
}

In all of our example, Menger-bounded groups will also be Scheepers-bounded,
and groups with the Menger property will also have the Scheepers property.

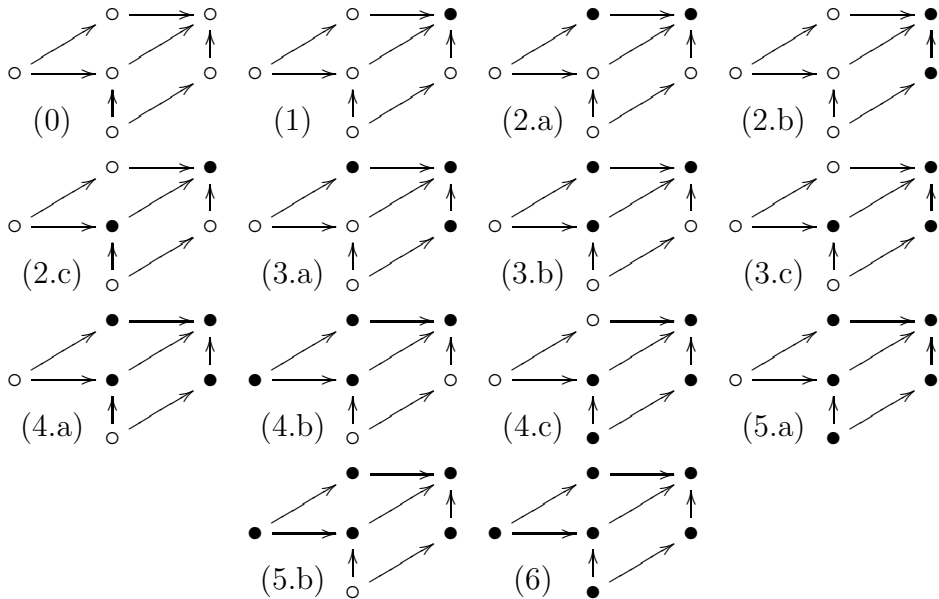
\begin{figure}[!ht]
$\xymatrix@R=10pt{
                      & \circ   \ar[r]        & \circ         \\
\circ   \ar[r]\ar[ur] & \circ   \ar[ur]       & \circ   \ar[u]\\
\hl{(0)}              & \circ   \ar[u]\ar[ur]
}$
$\xymatrix@R=10pt{
                      & \circ   \ar[r]        & \bullet       \\
\circ   \ar[r]\ar[ur] & \circ   \ar[ur]       & \circ   \ar[u]\\
\hl{(1)}                      & \circ   \ar[u]\ar[ur]
}$
$\xymatrix@R=10pt{
                      & \bullet \ar[r]        & \bullet       \\
\circ   \ar[r]\ar[ur] & \circ   \ar[ur]       & \circ   \ar[u]\\
\hl{(2.a)}                      & \circ   \ar[u]\ar[ur]
}$
$\xymatrix@R=10pt{
                      & \circ   \ar[r]        & \bullet       \\
\circ   \ar[r]\ar[ur] & \circ   \ar[ur]       & \bullet \ar[u]\\
\hl{(2.b)}            & \circ   \ar[u]\ar[ur]
}$
$\xymatrix@R=10pt{
                      & \circ   \ar[r]        & \bullet       \\
\circ   \ar[r]\ar[ur] & \bullet \ar[ur]       & \circ   \ar[u]\\
\hl{(2.c)}                      & \circ   \ar[u]\ar[ur]
}$
$\xymatrix@R=10pt{
                      & \bullet \ar[r]        & \bullet       \\
\circ   \ar[r]\ar[ur] & \circ   \ar[ur]       & \bullet \ar[u]\\
\hl{(3.a)}                      & \circ   \ar[u]\ar[ur]
}$
$\xymatrix@R=10pt{
                      & \bullet \ar[r]        & \bullet       \\
\circ   \ar[r]\ar[ur] & \bullet \ar[ur]       & \circ   \ar[u]\\
\hl{(3.b)}                      & \circ   \ar[u]\ar[ur]
}$
$\xymatrix@R=10pt{
                      & \circ   \ar[r]        & \bullet       \\
\circ   \ar[r]\ar[ur] & \bullet \ar[ur]       & \bullet   \ar[u]\\
\hl{(3.c)}                      & \circ   \ar[u]\ar[ur]
}$
$\xymatrix@R=10pt{
                      & \bullet \ar[r]        & \bullet       \\
\circ   \ar[r]\ar[ur] & \bullet \ar[ur]       & \bullet \ar[u]\\
\hl{(4.a)}                      & \circ   \ar[u]\ar[ur]
}$
$\xymatrix@R=10pt{
                      & \bullet \ar[r]        & \bullet       \\
\bullet \ar[r]\ar[ur] & \bullet \ar[ur]       & \circ   \ar[u]\\
\hl{(4.b)}                      & \circ   \ar[u]\ar[ur]
}$
$\xymatrix@R=10pt{
                      & \circ   \ar[r]        & \bullet       \\
\circ   \ar[r]\ar[ur] & \bullet \ar[ur]       & \bullet   \ar[u]\\
\hl{(4.c)}                      & \bullet \ar[u]\ar[ur]
}$
$\xymatrix@R=10pt{
                      & \bullet \ar[r]        & \bullet       \\
\circ   \ar[r]\ar[ur] & \bullet \ar[ur]       & \bullet   \ar[u]\\
\hl{(5.a)}            & \bullet \ar[u]\ar[ur]
}$
$\xymatrix@R=10pt{
                      & \bullet \ar[r]        & \bullet       \\
\bullet \ar[r]\ar[ur] & \bullet \ar[ur]       & \bullet \ar[u]\\
\hl{(5.b)}                      & \circ   \ar[u]\ar[ur]
}$
$\xymatrix@R=10pt{
                      & \bullet \ar[r]        & \bullet       \\
\bullet \ar[r]\ar[ur] & \bullet \ar[ur]       & \bullet \ar[u]\\
\hl{(6)}                      & \bullet \ar[u]\ar[ur]
}$
\caption{All settings consistent with the present arrows}\label{settings}
\end{figure}

To make our assertions visually clear,
we will use copies of the diagram with ``$\bullet$'' placed in positions
corresponding to a property that a given group has, and ``$\circ$'' in positions
corresponding to a property that the group does not have.
There are exactly $14$ possible settings to check. They all appear
in Figure \ref{settings}, where each setting is labelled $(n.m)$ where $n$ is the
number of $\bullet$'s in that setting.

Let us begin with the immediate examples.

\begin{prop}
Setting (0) is realized by $\Bgp$.
Setting (6) is realized by all countable subgroups of $\Bgp$.
\end{prop}
\begin{proof}
(0) See Theorem \ref{MenBdd}(4).

(6) All critical cardinalities are uncountable (Theorem \ref{crits}).
\end{proof}

A less trivial realization of Setting (6) is the group $\tilde G$ in the proof of
Theorem \ref{RBNotMen} constructed under the hypothesis $\fb\le\cov(\cM)$ (see Remark \ref{RBNotM2}).
Note further that any dominating subset of $\Bgp$ generates a realization of
the second setting (Theorem \ref{MenBdd}).

Section \ref{critcards} implies consistent realizations of Settings (4.b), (4.c), and (2.c):
For Setting (4.b), assume that $\cov(\cM)<\fb$ (which is consistent). By Theorem \ref{crits}(3),
there is a subgroup $G$ of $\Bgp$ such that $|G|=\cov(\cM)$ and $G$ is not Rothberger-bounded.
As $|G|<\fb$, we have by By Theorem \ref{crits}(2) that $G$ has the Hurewicz property.
Setting (4.c) follows, in a similar manner, from the consistency of $\fb<\cov(\cM)$.
For Setting (2.c), assume that $\max\{\fb,\cov(\cM)\}<\fd$ (which is consistent), and take
by Theorem \ref{crits}(2,3) subgroups $G,H$ of $\Bgp$ such that $|G|=\fb$, $|H|=\cov(\cM)$, and
$G$ is not Hurewicz-bounded and $H$ is not Rothberger-bounded.
Let $M=G+H$. $|M|=\max\{\fb,\cov(\cM)\}<\fd$, and thus by \ref{crits}(1) $M$ has
the Menger property. As $M$ contains $G$ and $H$, it is not Hurewicz-bounded nor Rothberger-bounded.

By Laver's result, each of the settings (2.b), (3.a), (3.c), (4.a), (4.c), (5.a), and (5.b)
requires some special hypothesis to be realized. Settings (4.b) and (2.c), mentioned in the
previous paragraph, are not in that list, and indeed do not require any special hypothesis.

\begin{thm}\label{MNotHBNotRB}
Settings (4.b) and (2.c) are realized.
\end{thm}
\begin{proof}
(4.b) The proof of Theorem \ref{cantorgroup} shows that
the group $G=\langle \Cantor\rangle$ is as required.

(2.c) By Theorem \ref{cantorgroup}, $G=\langle \Cantor\rangle$ is $\sigma$-compact.

For a proof of the following, see Lemma 2.6 of \cite{AddQuad}.

\begin{lem}[folklore]\label{sigmaPlus}
Assume that $X$ is a $\sigma$-compact space, and $Y$ has the Menger property.
Then $X\x Y$ has the Menger property.\hfill\qed
\end{lem}

The following also follows from Theorems \ref{MenBdd}, \ref{SchBdd}, and \ref{HurBdd}.

\begin{lem}[Babinkostova \cite{Bab05}]
Assume that $G$ is a Hurewicz-bounded subgroup of $\Bgp$, and $H$ is a Menger-
(respectively, Scheepers-, Hurewicz-) bounded
subgroup of $\Bgp$.
Then $G+H$ is Menger- (respectively, Scheepers-, Hurewicz-) bounded.\hfill\qed
\end{lem}

On the other hand, $G+H$, containing $G$, is not Rothberger-bounded.

Take $G=\langle \Cantor\rangle$ and let $H$ be as in Theorem \ref{CPthm}.
By Lemma \ref{sigmaPlus}, $G+H$ (a continuous image of $G\x H$) has the Menger property.
Containing $G$, it is not Rothberger-bounded, as is evident from Theorem \ref{RothBdd}(2).
Containing $H$, it is not Hurewicz-bounded, either.
\end{proof}

Recall that Setting (4.c) is realized when $\max\{\fb,\cov(\cM)\}<\fd$.
This setting can also be realized when the involved critical
cardinalities are equal.

\begin{thm}\label{RNotHB}
Assume that $\cov(\cM)=\fc$. Then Setting (4.c) is realized.
\end{thm}
\begin{proof}
See the proof of Theorem \ref{RBNotHB}.
\end{proof}

\begin{thm}\label{HBNotMNRB}
Setting (2.a) is realized.
\end{thm}
\begin{proof}
Let $G$ be as in Theorem \ref{HBNotM} and take $\langle \Cantor\rangle+G$.
Use Lemma \ref{sigmaPlus}.
\end{proof}

\begin{thm}
Assume that $\cov(\cM)=\fb=\fd$. Then Setting (3.a) is realized.
\end{thm}
\begin{proof}
See the proof of Theorem \ref{RBNotMen}.
\end{proof}

\begin{thm}
Setting (1) is realized.
\end{thm}
\begin{proof}
Take $G+H+\langle \Cantor\rangle$, where $G$ is from Theorem \ref{HBNotM}
and $H$ is from Theorem \ref{CPthm}.
As this group contains all three groups, it cannot satisfy any
of our properties, except perhaps Menger-boundedness.
As $H$ is Menger-bounded and $\langle \Cantor\rangle$ is $\sigma$-compact,
$H+\langle \Cantor\rangle$ is Menger-bounded. As $G$ is a subgroup of a $\sigma$-compact group,
$G+H+\langle \Cantor\rangle$ is Menger-bounded.
\end{proof}

\begin{thm}\label{Hamel1}
Setting (3.b) is realized.
\end{thm}
\begin{proof}
By a result of Mycielski \cite{Myc64},
there is a subset $C$ of $\Cgp$, which is linearly independent
over $\Z_2$, and is homeomorphic to $\Cgp$.
Then $C$ contains a subset homeomorphic to $\Bgp$.
Let $X$ be as in Theorem \ref{CPthm},
and let $Y$ be a topological embedding of $X$ in $C$.
Let $G=\langle Y\rangle$.
As all finite powers of $Y$ have the Menger property,
$G$ has this property too.

As $C$ is closed and $Y=G\cap C$, $Y$ is closed in $G$.
As $Y$ does not have the Hurewicz property,
$G$ does not have the Hurewicz property.

By Lemmas \ref{powers} and \ref{powers2}, the subgroup
$\tilde G=\langle G\rangle$ of $\Bgp$ has the Menger property in all finite
powers, and does not have the Hurewicz property.
It is Hurewicz-bounded since it is a subgroup of the $\sigma$-compact
group $\langle \Cantor\rangle$.
Take the group $\langle \Cantor\rangle+\tilde G$ and use Lemma \ref{sigmaPlus}.
\end{proof}

\begin{thm}\label{HamelLuzin}
Assume that $\cov(\cM)=\fc$. Then Setting (5.a) is realized.
\end{thm}
\begin{proof}
This is similar to the proof of Theorem \ref{Hamel1}.
Let $C$ be as in that proof, and $L$ be as in Theorem \ref{RBNotHB}.
Let $M$ be a topological embedding of $L$ in $C$, and take $G=\langle M\rangle$ in $\PN$.
As all finite powers of $M$ have the Rothberger property,
this holds for $G$ too (Lemma \ref{powers}).

As $M$ is a closed subset of $G$ and
$M$ does not have the Hurewicz property,
$G$ does not have the Hurewicz property.

Take $\tilde G=\langle G\rangle$ in $\Bgp$, and use Lemmas \ref{powers} and \ref{powers2}.
\end{proof}

\begin{thm}\label{<d}
Assume that $\fb=\cov(\cM)<\fd$.
Then Setting (4.a) is realized.
\end{thm}
\begin{proof}
We first get a subgroup of $\PN$.

\begin{lem}\label{bscale}
Assume that $\fb=\cov(\cM)$. Then there is a family
$B=\{b_\alpha : \alpha<\fb\}\sbst\roth$ such
that the subgroup $G=\langle B\rangle$ of $\PN$ is contained
in $\roth\cup\{\emptyset\}$, and $e[B]\sbst\NN$ has the following properties:
\be
\itm $e[B]=\{e(b_\alpha) : \alpha<\fb\}$ is $\le^*$-increasing
with $\alpha$ and $\le^*$-unboun\-ded,
\itm There is a continuous $\Phi:\NN\to\NN$ such that
$\Phi[e[B]]$ is a witness for the combinatorial definition of $\cov(\cM)$.
\ee
\end{lem}
\begin{proof}
Let $\kappa=\fb=\cov(\cM)$, and $Y=\{y_\alpha : \alpha<\kappa\}\sbst\NN$ be
such that there is no $f\in\NN$ such that for each $y\in Y$, $f(n)=y(n)$ for some $n$.

Fix a $\le^*$-unbounded subset $\{f_\alpha : \alpha < \kappa\}$ of $\NNup$.
Fix a bijection $\iota:\N\x\N\to\N$ such that $m,n\le\iota(m,n)$ for all $m,n$.\footnote{
The canonical bijection $\iota:\N\x\N\to\N$ has the required property.
}
We will use the homeomorphism $\Psi:\NN\x\NN\to\NN$ defined by
$$\Psi(f,g)(n)=\iota(f(n),g(n))$$
for each $n$.
Note that $h\le^*\Psi(f,g)$ whenever $h\le^* g$.

Define $B=\{b_\alpha : \alpha<\kappa\}\sbst\roth$ by induction on $\alpha$:
At step $\alpha$,
let
$$G_\alpha=\langle \{b_\beta : \beta<\alpha\}\cup\Fin\rangle = \langle \{b_\beta : \beta<\alpha\}\rangle\Delta\Fin.$$
Let $h_\alpha$ be a $\le^*$-bound of
$\{e[b_\beta] : \beta<\alpha\}$. We may assume that $f_\alpha\le h_\alpha$.
Let
$$A_\alpha = \{\Psi(y_\alpha, g) : g\in\NNup,\ h_\alpha\le^* g\}.$$
As $\Psi$ is injective, $|A_\alpha|=\fc\ge\kappa>|G_\alpha|$,
so we can choose $g_\alpha\in A_\alpha\sm e[G_\alpha\cap\roth]$.
Take $b_\alpha = e\inv(g_\alpha)=\im(g_\alpha)$.

Clearly, $\langle B\rangle\sbst\roth\cup\{\emptyset\}$ and (1) holds.
For (2), let $\Phi:\NN\to\NN$ be the mapping
$\Phi(\Psi(y,g)) = y$.
Then $\Phi[e[B]]=Y$.
\end{proof}

Let $\kappa=\fb=\cov(\cM)$.
Take $B\sbst\roth$ as in Lemma \ref{bscale}.
As $\fb\le\cov(\cM)$, all finite powers of $B\cup\Fin$ have the Rothberger
property \cite{ideals}, and by Lemma \ref{powers}, so does $G=\langle B\cup\Fin\rangle$ in $\PN$.
Consequently, its subgroup $H=\langle B\rangle$ is Rothberger-bounded.
As $|H|=\kappa<\fd$, $H$ has the Menger property.

As $e[H\sm\{\emptyset\}]\spst e[B]$ is $\le^*$-unbounded, $H\sm\{\emptyset\}$ does not have
the Hurewicz property. By Lemma \ref{removept}, $H$ does not have
the Hurewicz property.
Now, $\Phi[e[H\sm\{\emptyset\}]]\spst \Phi[e[B]]$,
and $\Phi[e[B]]$ does not have strong measure zero.
As having strong measure zero is a hereditary property,
$\Phi[e[H\sm\{\emptyset\}]]$ does not have strong measure zero.
As $\Phi$ is continuous, $e[H\sm\{\emptyset\}]$ does not have the Rothberger property,
and therefore neither does $H\sm\{\emptyset\}$. Consequently,
$H$ does not have the Rothberger property.

To get a subgroup of $\Bgp$ as required, take $\tilde G=\langle G\rangle$
and $\tilde H=\langle H\rangle$ in $\Bgp$.
By Lemma \ref{powers}, $\tilde G$ has the Rothberger property,
and therefore $\tilde H$ is Rothberger-bounded. Being a subgroup of $\langle \Cantor\rangle$,
it is  Hurewicz-bounded.
It has the Menger property because $|\tilde H|=\kappa<\fd$.
By Lemma \ref{powers2}, $\tilde H$ does not have the Hurewicz nor the Rothberger property.
\end{proof}

Having successfully realized all but three of the settings in Figure \ref{settings},
one may be tempted to assume that all settings can be realized.
Surprisingly, this is not entirely the case.

\begin{thm}\label{NSW}
There does not exists a realization of Setting (5.b).
\end{thm}
\begin{proof}
In \cite{NSW} it is proved that if a set of reals has the Hurewicz
property and has strong measure zero, then it has the Rothberger
property (see \cite{prods} for a simple proof of that assertion).
Use Theorem \ref{RothBdd}.
\end{proof}

Only the Settings (2.b) and (3.c) of Figure \ref{settings} remain unsettled.

\begin{prob}
Is it consistent that there is a subgroup $G$ of $\Bgp$ such that
$G$ has strong measure zero, is unbounded (with respect to $\le^*$), and
does not have the Menger property?
\end{prob}

\begin{prob}
Is it consistent that there is a subgroup $G$ of $\Bgp$ such that
$G$ has strong measure zero and Menger's property, but is unbounded (with respect to $\le^*$) and
does not have the Rothberger property?
\end{prob}

\subsection*{Acknowledgments}
We thank Heike Mildenberger, Lyubomyr Zdomskyy, and the referee for their useful comments.
We also thank Taras Banakh, Lyubomyr Zodmskyy, and Petr Simon,
for suggesting the elimination of one quantifier in Theorems \ref{MenBdd} and \ref{SchBdd}.

The major part of the research towards this paper was carried out
when the first author visited Bar-Ilan University, and the second
author was a Koshland Scholar at the Weizmann Institute of Science.
The authors thank the Mathematics Faculties of both institutions,
and especially their hosts Boris Kunyavski, Adi Shamir, and Gideon Schechtman,
for their very kind hospitality.

\ed